\input amstex
\documentstyle{amsppt}
\magnification=1200
\loadbold

\def\bdot{\boldsymbol\cdot}
\def\bm{\boldsymbol \mu}
\def\bsq{\boxtimes}
\def\t{\theta}
\def\h{{\sssize 1\!/\!2}}
\def\qh{q^{\h}}
\def\qm{q^{-\h}}
\def\th{t^{\h}}
\def\Z{\Bbb Z}
\def\N{\Bbb N}
\def\C{\Bbb C}

\def\l{\lambda}
\def\xs{x^*}
\def\ys{y^*}
\def\1{{\pm 1}}
\def\L{\Lambda}
\def\Lt{\L_{t,s}}
\def\Lp#1{\L_{t,s}^{\sssize #1}}
\def\Lh#1{\L_{t,s\th}^{\sssize #1}}
\def\k{\Bbbk}
\def\kx{\k[x_1^\1,\dots,x_n^\1]}
\def\P{P^*}
\def\Pm{P^*_\mu}
\def\co{\operatorname{const}}
\def\d{\,d_q}
\def\dd#1{\,\frac{d_q{#1}}{#1}}
\def\inv{^{-1}}
\def\de{\delta}
\def\lr{\leftrightarrow}
\def\b{\beta_{q,t}}
\def\tht{\thetag}
\def\I{I}

\def\hV{\widehat{V}}
\def\hPi{\widehat{\Pi}}
\def\hD{\widehat{D}}
\def\hC{\widehat{C}}

\def\hx{\hat x}
\def\hy{\hat y}
\def\prs{\curlyeqprec}
\def\st{\,|\,}

\def\sq{\square}

\def\i{{}^{\sssize (@!@!@!@!@! i @!@!@!@!@!)}@!@!@!@!@!@!@!}
\def\ii#1{{}^{\sssize (@!@!@!@!@! #1 @!@!@!@!@!)}@!@!@!@!@!@!@!}
\def\D{\frak D}
\def\Dl{{\Cal D}_l}
\def\T{{\sssize T}}
\def\pmu{{}^\backprime\!\mu}
\def\Ugln{{\Cal U}(\frak{gl}(n))}
\def\UglN{{\Cal U}(\frak{gl}(N))}
\def\bc{\sssize \blacksquare}
\def\wc{\sssize \square}
\def\p{\psi}
\def\tm{{\widetilde{\mu}}}
\def\teta{{\widetilde{\eta}}}

\leftheadtext{Andrei Okounkov}
\rightheadtext{Interpolation Macdonald polynomials}

\topmatter

\title
$\boldkey B\boldkey C$-type interpolation Macdonald polynomials 
and binomial formula for Koornwinder polynomials
\endtitle

\author
Andrei Okounkov
\endauthor
\address
Department of Mathematics, University of Chicago,  5734 
University Avenue, Chicago, IL 60637-1546
\endaddress
\email
okounkov\@math.uchicago.edu
\endemail

\abstract
We consider 3-parametric polynomials $\Pm(x;q,t,s)$
which replace the $A_n$-series interpolation Macdonald polynomials $\Pm(x;q,t)$
for the $BC_n$-type root system. For these polynomials 
we prove an integral representation, a combinatorial
formula, Pieri rules, Cauchy identity, and we also show
that they do not satisfy any rational $q$-difference
equation. As $s\to\infty$ the polynomials $\Pm(x;q,t,s)$
become $\Pm(x;q,t)$. We also prove a binomial formula for
6-parametric Koornwinder polynomials.
\endabstract

\toc
\widestnumber\head{A2}
\head 0. Introduction \endhead 
\head 1. Definition and normalization \endhead
\head 2. Elementary properties \endhead
\head 3. Integral representation \endhead
\head 4. Highest degree term \endhead
\head 5. Combinatorial formula \endhead
\head 6. Pieri rules and Cauchy identity \endhead
\head 7. Binomial formula for Koornwinder polynomials \endhead
\head A1. Appendix 1. $q$-integrals \endhead
\head A2. Appendix 2. Absence of rational $q$-difference equations \endhead
\endtoc

\endtopmatter

\head
0.~Introduction
\endhead

The higher Capelli identities \cite{N,Ok1-2} 
describe two faces of a distinguished linear basis in the center 
of $\Ugln$. The elements of this basis, which are called {\it quantum
immanants},  have also a third and a very important face. It is
their Harish-Chandra image called {\it shifted Schur functions}, 
or $s^*\!$--functions for short.  The $s^*\!$--functions have 
remarkably many applications, see \cite{OO} and, for
example, also \cite{KOO}. In
particular, they were used in the proof of the higher Capelli
identities given in \cite{Ok1}.

Since the definition of the $s^*\!$--functions involves only
the Weyl group and the half-sum of positive roots $\rho$,
they can be generalized in a quite far reaching way. This was done
in the $A_n$ case 
by F.~Knop,  G.~Olshanski, S.~Sahi, and the author in 
\cite{KS,Kn,S,Ok3,OO3} and other papers, see References.
This theory  is remarkably parallel to the theory of ordinary Macdonald
polynomials $P_\mu(x;q,t)$ for the root system $A_n$.
In some aspects, such as e.g.\ duality \cite{Ok3}, \S6,  or binomial formula
\cite{Ok4}, it looks even more natural. 

The definition of the polynomials $\Pm(x;q,t)$, which generalize
shifted Schur functions in the way that
$$
s^*_\mu(x)=\lim_{q\to 1} \frac{\Pm(q^x;q,q)}{(q-1)^{|\mu|}}\,, 
$$
where
$$
q^x=(q^{x_1},\dots,q^{x_n})
$$
and $|\mu|$ stands for the number of squares in $\mu$, 
is the following. By definition, the polynomial $\Pm(x;q,t)$ is
the unique, up to a scalar, polynomial of degree $|\mu|$
that is symmetric in variables
$$
x_1 t^{n-1}, x_2 t^{n-2}, \dots, x_n\,, \tag 0.1
$$
and vanishes at the points
$$
\Pm(q^\l;q,t)=0\,, \quad \mu\not\subset\l\,,
$$
where $\l=(\l_1,\dots,\l_n)$ is a partition. 
We call these
polynomials {\it interpolation Macdonald polynomials}. They
are also known as shifted or inhomogeneous Macdonald polynomials. 

The  polynomials $\Pm(x;q,t)$ are the Harish-Chandra
image
$$
\D @>\quad\text{Harish-Chandra isomorphism }\quad>> 
\C(q,t)[x_1 t^{n-1}, x_2 t^{n-2}, \dots, x_n]^{S(n)}
$$
of some distinguished linear basis of the 
commutative algebra
$\D$ of Macdonald $q$-difference operators. 
Recall that the algebra $\D$ acts diagonally in the basis of Macdonald
polynomials
and the Harish-Chandra homomorphism associates to any
element of $\D$ its eigenvalues on the Macdonald
polynomials $P_\mu(x;q,t)$ viewed as a function
of $q^\mu$.

Informally, by analogy with \cite{Ok1},
one can think of the polynomials $\Pm(x;q,t)$
as of the Harish-Chandra image of certain Capelli-type
Laplace operators on  on some  fictitious 
``quantum symmetric space'' with an $A_n$-type
restricted root system. In the three cases 
$$
t=q^\t\,, \quad \t=\frac12,1,2
$$
there exist certain quantum symmetric spaces (see \cite{No2})
which are $q$-analogs of 
$$
\alignat2
\vspace{-3 \jot}
&U(n)/O(n)\,, &\qquad &\t=1/2\,,\\
&U(n)\,,         &\qquad &\t=1\,,\\
&U(2n)/Sp(2n)\,, &\qquad &\t=2\,.
\endalignat
$$

Remark, that from the definition it is not obvious at all
that the polynomials $\Pm(x;q,t)$ are anyhow related to ordinary Macdonald
polynomials. There is, therefore, a natural
question if any such relation
exists for other root systems. The root system
of maximal interest for applications is the non-reduced
root system of type $BC_n$, where we have even more
``live'' symmetric spaces than in the $A_n$ case. 
In \cite{OO2}, we considered the analogs of shifted
Schur functions for other classical groups. 

In this paper we study the following natural $BC_n$-type
analogs of the interpolation Macdonald polynomials $\Pm(x;q,t)$. 
By definition, the polynomial
$$
\Pm(x;q,t,s)\in\C(q,t,s)[x_1^\1,\dots,x_n^\1]
$$
is the unique, up to a scalar, polynomial of degree
$|\mu|$ which is symmetric in variables \tht{0.1} and also 
invariant under the transformations
$$
x_i t^{n-i} s \mapsto \frac1{x_i t^{n-i} s}\,,
$$
such that
$$
\Pm(q^\l;q,t,s)=0 \quad \text{if} \quad \mu\not\subset\l 
$$
for any partition $\l$. The vector 
$\rho=(\rho_1,\dots,\rho_n)$ such that
$$
q^\rho=(t^{n-1} s,\dots,t s,s)
$$
plays the role of the half-sum of positive roots
taken with multiplicities. 

Our analysis shows that this  $BC_n$-type case 
has some new and unexpected features.

We prove that the polynomials $\Pm(x;q,t,s)$ are
actually  very close to the polynomials $\Pm(x;q,t)$.
In particular, as $s\to\infty$ all negative powers
of variables disappear and
$$
\Pm(x;q,t,s)\to\Pm(x;q,t), \quad s\to\infty\,. \tag 0.2
$$
Moreover, the new parameter $s$ can be inserted into the two
explicit formulas for $\Pm(x;q,t)$ obtained in \cite{Ok3}, namely, the {\it
integral representation\/} and the {\it combinatorial
formula}, in such a way that \tht{0.2}
becomes evident. This is done in sections 3--5, where the proofs
are suitable modifications
of the proofs from \cite{Ok3}. There is, however, also a major
difference between the polynomials $\Pm(x;q,t,s)$ and
$\Pm(x;q,t)$, see below.

Recall  that in the 
$s^*\!$--functions case the combinatorial formula
and the integral representation, which was called in that
context the {\it coherence property}, have a clear
interpretation in terms of the higher Capelli 
identities. The combinatorial formula reflects 
the formula \tht{3.24} in \cite{Ok1} for the 
LHS of the higher Capelli identities, see \cite{Ok1}, \S3.7.
The coherence property describes the action 
of the averaging operator
$$
\text{center of $\Ugln$}\quad
@>\quad\text{averaging}\quad>>
\quad \text{center of $\UglN$}\,, \quad n<N 
$$
in the basis of quantum immanants. This action
can be easily computed either by explicit
averaging of the differential operator in the
RHS of the higher Capelli identities, see \cite{Ok1}, \S5.1, 
or by using the vanishing property, see \cite{OO}, \S10. 

The  integral representation \tht{3.1} 
is in fact a natural analog of
the classical Weyl character formula. The denominator $D(x^*;q,t)$ in
\tht{3.1} is a straightforward generalization of the 
product in the denominator of the Weyl formula.  The
alternating sum in the numerator of the Weyl formula
turns into  a multivariate $q$-integral which is, by 
definition, alternating sum of contributions of
lower and upper bounds of integration. One proves that
the numerator is divisible by the denominator
using, essentially, the fact that any $q$-integral
$$
\int_v^u z^l \,d_q z\,, \quad l\ne -1\,,
$$
is divisible by $u-v$ in $\C[u^\1,v^\1]$, see Appendix 1.

On the other hand, the relation of the polynomials 
$\Pm(x;q,t,s)$ to the Macdonald polynomials for the 
$BC_n$ root system is only an indirect one, via the 
{\it binomial theorem}. 
In fact, it is natural to work with {\it Koornwinder polynomials}
\cite{K}, which generalize Macdonald polynomials for classical root
systems and have as many as $6$ parameters. We show that using
the polynomials $\Pm(x;q,t,s)$ one obtains a binomial
formula for Koornwinder polynomials just along the same lines
as in \cite{Ok4}. The only property of the polynomials $\Pm(x;q,t,s)$ needed
in the proof of the binomial formula, as an abstract identity
of symmetric polynomials, is the identification \tht{4.1} of the
top homogeneous term of $\Pm(x;q,t,s)$.

As an application of this binomial formula, 
or, more precisely, of its classical limit, one
obtains a BC analog of the results of \cite{OO4}.
The details will be given in a forthcoming paper by
G.~Olshanski and the author. 

We conclude this introduction with remarks on some
open problems and also some {\it no go} results.  

By the binomial theorem, see \cite{OO3}, \S6, and \cite{Ok4}, 
the  integral representation of $\Pm(x;q,t)$ reflects the {\it
branching rule} for the ordinary $A_n$-type Macdonald polynomials.
It is likely that the integral representation 
of the polynomials $\Pm(x;q,t,s)$
reflects branching rules for the Koornwinder polynomials.
The 
$$
s\mapsto s \th
$$
shift in the integral representation \tht{3.1} suggests
that the branching rules for  Koornwinder polynomials
must have two natural steps just like the branching rules for
characters of orthogonal and symplectic series.

The combinatorial formula for the $A_n$-type 
polynomials $\Pm(x;q,t)$ and the integral representation
of $\Pm(x;q,t)$ are equivalent to each other by
virtue of the following symmetry, see \cite{Ok3}, 
formula \tht{2.1}, 
$$
\frac{\Pm(a q^{-\l_n}, \dots, a q^{-\l_1}; q,t)} 
{\Pm(a, \dots, a ; q,t)} =
\frac{P^*_\l(a q^{-\mu_n}, \dots, a q^{-\mu_1}; q,t)} 
{P^*_\l(a, \dots, a ; q,t)} \tag 0.3
$$
between the argument $x$ and the label $\mu$ of $\Pm(x;q,t)$.
Here $a$ is an arbitrary number. As $a\to\infty$ the symmetry
\tht{0.3} becomes the symmetry for ordinary Macdonald polynomials.

It looks like no such  symmetry connects the
integral representation \tht{3.1} to  the  combinatorial 
formula  \tht{5.3} for the $BC_n$-type interpolation polynomials.
First of all, the symmetry \tht{0.3} uses the automorphism of the weight lattice
$$
\align
\l&\mapsto - w_0(\l)\,,\\
(\l_1,\dots,\l_n) &\mapsto (-\l_n,\dots,-\l_1)\,,
\endalign
$$
which is trivial in the $BC_n$ case. Here $w_0$ stands for
the longest element of the Weyl group. Another obstruction to symmetry 
is that the weight function $\p_{\l/\mu}$
in the combinatorial formula, which is the same $\p_{\l/\mu}$
as in  \cite{M}, VI.6.24 and Example VI.6.2.(b), 
$$
\multline
\p_{\l/\mu}=
\prod_{i\le j\le\ell(\mu)}
\frac
{(q^{\mu_i-\mu_j} t^{j-i+1})_\infty}
{(q^{\mu_i-\mu_j+1} t^{j-i})_\infty}
\frac
{(q^{\l_i-\l_{j+1}} t^{j-i+1})_\infty}
{(q^{\l_i-\l_{j+1}+1} t^{j-i})_\infty} \times \\
\frac
{(q^{\l_i-\mu_j+1} t^{j-i})_\infty}
{(q^{\l_i-\mu_j} t^{j-i+1})_\infty}
\frac
{(q^{\mu_i-\l_{j+1}+1} t^{j-i})_\infty}
{(q^{\mu_i-\l_{j+1}} t^{j-i+1})_\infty}
\endmultline \tag 0.4
$$
is, unlike the weight function in the integral representation
\tht{3.1}, not invariant under transformations of the form
$$
q^{\l_i} t^{n-i} s \mapsto 
\frac 1{q^{\l_i} t^{n-i} s} \,.
$$

This problem is closely related to the following 
important difference between the 
polynomials $A_n$ and $BC_n$ case.
The polynomials $\Pm(x;q,t,s)$
{\it do not} satisfy any $q$-difference
equations with rational coefficients, see Appendix 2. 
Recall that there are very useful $q$-difference equations
for $\Pm(x;q,t)$ discovered in \cite{K,S2}. Explicit formulas
for higher order difference equations can be found in \cite{Ok4},
\S3. It seems to be a very interesting and important problem to
find a geometric or representation-theoretic meaning of the Knop-Sahi
difference equation and to understand why they do not have
obvious $BC$-analogs. 

Finally, it is natural to ask if one can add some more
parameters to the 3-parametric polynomials $\Pm(x;q,t,s)$
and still preserve some of their properties. It has been
shown in \cite{Ok5} that the property
$$
\Pm(q^\l;q,t,s)=0\,, \quad \mu\not\subset\l\,,
$$
which requires the polynomials $\Pm(x;q,t,s)$ to
vanish at an infinite set of points $q^\l$, $\mu\not\subset\l$,
 characterizes
the polynomials $\Pm(x;q,t,s)$ inside a very general
class of interpolation polynomials. It follows from
this characterization, see \cite{Ok5}, that neither
integral representation nor combinatorial formula
exist for any more general interpolation polynomials
in that class. 

\head  
1.~Definition and normalization
\endhead

Let $q,t,s$ be three parameters. 
We shall assume that
$$
q^i t^j s^k \ne 1 \quad 
\text{for all}\quad i,j,k\in\Z_+\,.
$$
Set
$$
\k=\C(q,t,s)\,.
$$
We shall assume that
$$
q,t\in\C\quad\text{and}\quad |q|,|t|<1
$$
wherever convergence is involved.
We shall also need the square roots $\qh$ and $\th$.
Consider the $BC$-type Weyl group
$$
W=S(n)\ltimes \Z_2^n \,.
$$
and its standard action on $\kx$. Denote by 
$$
\Lt\subset\kx
$$
the subalgebra of polynomials that are $W$-invariant in new
variables
$$
\xs_i=x_i t^{n-i} s\,, \quad i=1,\dots,n\,. \tag 1.1
$$
Note that for any $f\in\Lt$
$$
\deg f \ge 0
$$
and $\deg f=0$ only for constant polynomials. Here $\deg f$
denotes the usual degree, that is the maximal total degree in
$x_1,\dots,x_n$ of all
monomials of $f$.

\definition{Definition 1.1}
Let $\mu$ be a partition with at most $n$ parts.
By definition, 
$$
\Pm(x_1,\dots,x_n;q,t,s) 
$$
is the element of $\Lt$ satisfying the following conditions:
\roster
\item $\deg \Pm(x;q,t,s) \le|\mu|$,
\item $\Pm(q^\l;q,t,s)=0$ if $\mu\not\subset\l$,
\item $\Pm(q^\mu;q,t,s)=H(\mu,n;q,t,s)$\,.
\endroster
\enddefinition

Here $H(\mu,n;q,t,s)$ is just a nonzero normalization constant
which we shall specify below.

It is clear that existence of $\Pm$ implies uniqueness.
Below we shall give two explicit formulas for $\Pm$.
It is also clear that

\proclaim{Proposition 1.1} The polynomials $\Pm(x;q,t,s)$, where
$\mu$ ranges over partitions with at most $n$ parts, form
a $\k$-basis of the vector space $\Lt$. The degree 
\roster
\item"(1${}'$)"  $\deg \Pm(x;q,t,s) = |\mu|$
\endroster
is exactly $|\mu|$.
\endproclaim

For any polynomial
$$
f(x)\in\Lt
$$
the coefficients $f_\mu$ in the expansion
$$
f(x)=\sum_\mu f_\mu \Pm(x;q,t,s)
$$
can be found from the following non-degenerate
triangular (with respect to the ordering of
partitions by inclusion) system of linear
equations
$$
f(q^\l)=\sum_\mu f_\mu \Pm(q^\l;q,t,s)\,.
$$
Here $\l$ ranges over all partitions with at most $n$ parts.

Now we shall specify the constant $H(\mu,n;q,t,s)$.
Consider the diagram
$$
\mu=(\mu_1,\dots,\mu_n)
$$
as the following skew diagram
$$
\mu=M/\bm\,,
$$
where 
$$
\align
M&=
(2\mu_1,\mu_1+\mu_2,\dots,\mu_1+\mu_n,\mu_1-\mu_n,\dots,\mu_1-\mu_2,0)\,,\\
\bm&=
(\mu_1,\mu_1,\dots,\mu_1,\mu_1-\mu_n,\dots,\mu_1-\mu_2,0)\,.
\endalign
$$
For the diagram
$$
\mu=(4,2,1,1)
$$
the diagram $M$ looks as follows (the black squares correspond to
the subdiagram $\bm$)
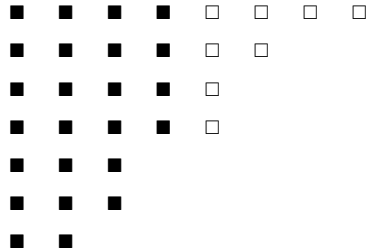
$$
\matrix
\bc&\bc&\bc&\bc&\wc&\wc&\wc&\wc\\
\bc&\bc&\bc&\bc&\wc&\wc\\
\bc&\bc&\bc&\bc&\wc\\
\bc&\bc&\bc&\bc&\wc\\
\bc&\bc&\bc\\
\bc&\bc&\bc\\
\bc&\bc
\endmatrix
$$
Each square
$$
\sq\in\mu
$$
has its mirror image $\bsq$ in the diagram $\bm$. If $\sq$ is
the $i$-th row and $j$-th column of $\mu$ then $\bsq$ is in the
$i$-th row and $(\mu_i-j+1)$-st column of $\bm$. The mirror images
of four squares $\sq\in\mu=(4,2,1,1)$ are given in the following 
picture
$$
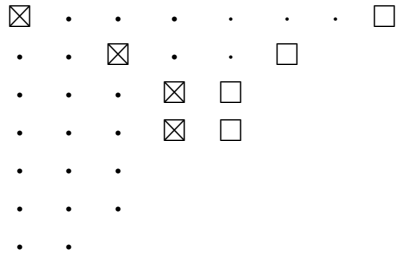
\matrix
\bsq &\bdot&\bdot&\bdot&\cdot&\cdot&\cdot&\square\\
\bdot&\bdot&\bsq&\bdot&\cdot&\square\\
\bdot&\bdot&\bdot&\bsq&\square\\
\bdot&\bdot&\bdot&\bsq&\square\\
\bdot&\bdot&\bdot\\
\bdot&\bdot&\bdot\\
\bdot&\bdot
\endmatrix
$$

Recall the following notations of Macdonald, see \cite{M}, Ch.~I. 
Given a partition $\mu$ set
$$
n(\mu)=\sum_i (i-1) \mu_i = \sum_j \mu'_j(\mu'_j-1)/2 \,.
$$
Recall that for each square $\sq=(i,j)\in\mu$ the numbers
$$
\alignat2
&a(\sq)=\mu_i-j,&\qquad &a'(\sq)=j-1,\\
&l(\sq)=\mu'_j-i,&\qquad &l'(\sq)=i-1,
\endalignat
$$
are called arm-length, arm-colength, leg-length, and
leg-colength respectively. The numbers
$$
a(\bsq)=\mu_i+a'(\sq)\,,\qquad l(\bsq)=l(\sq)+2(n-\mu'_j)
$$
are the arm-length and leg-length of the mirror image $\bsq$
of $\sq$
measured with respect to the big diagram $M$.

\definition{Definition 1.2} $H(\mu,n;q,t,s)$ is the following
normalization constant
$$
\frac{t^{n(\mu)-2(n-1)|\mu|}}{q^{2n(\mu')+|\mu|}s^{2|\mu|}}
\prod_{\sq\in\mu}
\left(q^{a(\sq)+1} t^{l(\sq)} - 1\right)
\left(s^2 q^{a(\bsq)} 
t^{l(\bsq)}-1\right)\,.
$$
\enddefinition

Observe that 
$$
\align
\sum_{\sq\in\mu} a(\bsq) &= |\mu|+ 3n(\mu')\,,\\
\sum_{\sq\in\mu} l(\bsq) &= 2(n-1)|\mu|- 3n(\mu)\,.
\endalign
$$
Therefore,
$$
H(\mu,n;q,t,s)\to \frac{q^{n(\mu')}}{t^{2n(\mu)}}
\prod_{\sq\in\mu}
\left(q^{a(\sq)+1} t^{l(\sq)} - 1\right), \quad s\to\infty\,,
$$
which is the normalization constant for two-parameter interpolation
Macdonald polynomials, see \cite{Ok3}.

Recall also that
$$
\sum_{\sq\in\mu} a(\sq)=n(\mu')\,, \quad
\sum_{\sq\in\mu} l(\sq)=n(\mu)\,.
$$

Remark that $\Pm(x;q,t,s)$ depends only on $s^2$, not $s$, because
the polynomial $\Pm(x;q,t,-s)$ satisfies all conditions of the definition of
$\Pm(x;q,t,s)$.

\head 
2.~Elementary properties
\endhead

We have the three following elementary propositions.

\proclaim{Proposition 2.1}
Suppose $\mu_n>0$ and put 
$$
\mu-\bar1=(\mu_1-1,\dots,\mu_n-1)\,.
$$
Then
$$
\multline
\Pm(x;q,t,s)=\\q^{|\mu-\bar1|}  \prod_i (x_i t^{1-i} - t^{1-n})
\left(1-\frac1{x_i t^{n-i} s^2}\right) \,
\P_{\mu-\bar1}(x/q;q,t,s q)\,. 
\endmultline \tag 2.1
$$
\endproclaim

\demo{Proof}
It is clear the RHS of \tht{2.1} satisfies all conditions
of the definition 1.1 except for normalization. 
 Evaluate it  at $x=q^\mu$. We obtain
$$
\multline
q^{-|\mu|+|\mu-\bar1|} t^{-3n(n-1)/2} s^{-2n} \times \\
\prod_i \left(q^{\mu_i} t^{n-i} - 1\right)
\left(s^2 t^{n-i} q^{\mu_i}-1\right)\,
H(\mu-\bar1,n;q,t,s q) \,,
\endmultline 
$$
which equals $H(\mu,m;q,t,s)$ because
$$
\alignat2
&n(\mu)&&=n(\mu-\bar1)+n(n-1)/2\,,\\
&n(\mu')&&=n((\mu-\bar1)')+|\mu-\bar1|\,,\\
&|\mu|&&=|\mu-\bar1|+n\,.
\endalignat
$$
This concludes the proof.  \qed
\enddemo

A similar computation proves the following

\proclaim{Proposition 2.2} Suppose $\mu_n=0$. Then
$$
\Pm(x_1,\dots,x_{n-1},1;q,t,s)=
\Pm(x_1,\dots,x_{n-1};q,t,s t) \,. \tag 2.2
$$
\endproclaim

\proclaim{Proposition 2.3} 
$$
\Pm(1/x_1,\dots,1/x_{n};1/q,1/t,1/s)=
s^{2|\mu|} t^{(2n-2)|\mu|} \Pm(x_1,\dots,x_{n};q,t,s) \,.
$$
\endproclaim
\demo{Proof} 
Again, it is clear that the LHS satisfies all conditions
of the definition of $\Pm(x;q,t,s)$ except for normalization.
Compute $H(\mu,n;1/q,1/t,1/s)$. We obtain
$$
\frac{q^{2n(\mu')+|\mu|}s^{2|\mu|}}{t^{n(\mu)-2(n-1)|\mu|}}
\prod_{\sq\in\mu}
\left(q^{-a(\sq)-1} t^{-l(\sq)} - 1\right)
\left(s^{-2} q^{-a(\bsq)} 
t^{-l(\bsq)}-1\right)\,,
$$
which equals 
$$
\multline
\frac{t^{n(\mu)}}{q^{|\mu|+2n(\mu')}}
\prod_{\sq\in\mu}
\left(q^{a(\sq)+1} t^{l(\sq)} - 1\right)
\left(s^2 q^{a(\bsq)} 
t^{l(\bsq)}-1\right)
=\\
=s^{2|\mu|} t^{(2n-2)|\mu|} H(\mu,n;q,t,s)\,. \qed
\endmultline
$$
\enddemo

\head
3.~Integral representation
\endhead

In this section we obtain a $q$-integral representation of
$\Pm(x;q,t,s)$, which is a way to obtain the polynomial
$$
\Pm(x_1,\dots,x_n;q,t,s)
$$
in $n$ variables from the polynomial
$$
\Pm(y_1,\dots,y_{n-1};q,t,s\th)
$$
in $n-1$ variable. 

Basic facts about $q$-integrals are recalled in the Appendix.
Introduce some notations; set
$$
(a)_\infty=(1-a)(1-q a)(1-q^2 a)\dots \,.
$$
Consider the following products
$$
\align
&V(x_1,\dots,x_n)=\det 
\left[\left(x_i^{n-j+1} - x_i^{-(n-j+1)}
\right)\right]_{1\le i,j \le n}\,, 
\\
&\Pi(x_1,\dots,x_n;y_1,\dots,y_{n-1};q,t) =
\prod_{i,j} \frac {(\qh (x_i)^\1 (y_j)^\1)_\infty}
{(\th (x_i)^\1 (y_j)^\1)_\infty}\,,
\\
&D(x_1,\dots,x_n;q,t)=\prod_{i<j} ((x_i+x^{-1}_i)-
(x_j+x^{-1}_j))
\frac {(q (x_i)^\1 (x_j)^\1)_\infty}
{(t  (x_i)^\1 (x_j)^\1)_\infty}\,.
\endalign
$$
Set also
$$
\ys_i=y_i t^{n-i-\h} s, \quad i=1,\dots,n-1\,.
$$

We shall integrate over the domain
$$
\int_{y\prs x} \d y =
\int^{qx_1}_{x_2} \d y_1 \dots  
\int^{qx_{n-1}}_{x_{n}} \d y_{n-1}
$$
with respect to the following beta-type 
measure
$$
\b(d_qy)= V(\ys;q,t) \,\Pi(\xs;\ys;q,t) \, \dd{y_1} \dots \dd{y_{n-1}}\,.
$$
Here the variables $\xs$ are defined in \tht{1.1}. 
We have

\proclaim{Theorem 3.1 (Integral representation)} Suppose $\mu_n=0$. Then
$$
\Pm(x;q,t,s)=\frac{1}
{C(\mu,n)}\frac{t^{-|\mu|}}{D(\xs;q,t)} 
\int_{y\prs x} \b(d_qy) \, \Pm(y; q,t,s\th)\,. \tag 3.1 
$$
Here 
$$
\align
C(\mu,n)&=q^{n(n-1)/4} \prod_{i<n} B_q(\mu_i+(n-i)\theta,\theta)\\
&=q^{n(n-1)/4} (1-q)^{(n-1)} \prod_{i<n}
\frac{(q)_\infty}{(t)_\infty}
\frac{(q^{\mu_i} t^{n-i+1})_\infty}
{(q^{\mu_i} t^{n-i})_\infty}\,,
\endalign
$$
where $t=q^\theta$ and $B_q$ is the $q$-beta function \tht{A.3}. 
\endproclaim

This theorem together with proposition 2.1 gives a
recursive formula for all 3-parametric interpolation Macdonald
polynomials.

\proclaim{Proposition 3.2} The RHS of \tht{3.1}
is an element of $\Lt$ of degree $\le|\mu|$.
\endproclaim

In the proof it will be 
convenient to consider the following analogs
of the algebra $\Lt$. By
$$
\Lp{\pm\pm}\subset\kx
$$
denote the subalgebra of polynomials that are (anti)-invariant
with respect to permutation of $\xs_i$ and (anti)-invariant
with respect to transformations
$$
\xs_i\lr(\xs_i)\inv \,.
$$
In particular, $\Lp{++}=\Lt$. We shall also write
$$
\Lp{\pm\pm}[x]
$$
to stress the dependence on variables $x_1,\dots,x_n$.

\demo{Proof of proposition 3.2}
By analytic continuation we can assume that
$$
t=q^\theta, \quad \theta=2k+1, \quad k=0,1,2,\dots\,.
$$
In this case $\Pi(x;y;q,t)$ and $D(x;q,t)$ belong to  $\kx$.

Note that 
$$
\Pi(\xs;\ys;q,t)=0
$$
if 
$$
y_i=qx_i,q^2x_i,\dots,t q\inv x_i
$$
or if 
$$
y_i=q\inv x_{i+1},q^{-2}x_{i+1},\dots, q t\inv  x_{i+1} \,.
$$
Therefore, for any $i$, we can change the limits of integration
in \tht{3.1} as follows
$$
\int_{x_{i+1}}^{q x_i} \d y_i \,\bigg(\dots\bigg)=
\int_{q^{-s}x_{i+1}}^{q^{r+1} x_i} \d y_i \, \bigg(\dots\bigg)\,,
\quad r,s=1,\dots,2k\,. \tag  3.2
$$
In particular, we can replace our integration
by the following
$$
\int_{y\prs x}\,\bigg(\dots\bigg)=
\int^{(qt)^\h x_{1}}_{(q/t)^\h x_{2}} \dots
\int^{(qt)^\h  x_{n-1}}_{(q/t)^\h x_{n}}\,\bigg(\dots\bigg) \,.
$$

Denote by $f(x,y)$ the polynomial 
$$
f(x,y)=V(\ys;q,t) \,\Pi(\xs;\ys;q,t) \,\Pm(y; q,t,s\th)\,.
$$
We have
$$
f(x,y)\in\Lt[x], \quad f(x,y)\in\Lh{--}[y]\,.
$$
In other words, $f(x,y)$ is a $W$-antiinvariant polynomial
in variables $\ys_i$ with coefficients in $\Lt[x]$.

The following determinants $\Dl$ form a linear basis
in $\Lh{--}[y]$
$$
\Dl(y)=\det\left(
\matrix
&\vdots&\\
\hdots &
(\ys_i)^{l_j} - (\ys_i)^{-l_j}&
\hdots\\
&\vdots&
\endmatrix
\right)_{1\le i,j \le n-1}\,,
$$
where 
$$
l_1>l_2>\dots>l_{n_1}>0\,.
$$
Using \tht{A.2} we can 
evaluate the integral
$$
\int^{(qt)^\h x_{1}}_{(q/t)^\h x_{2}} \dots
\int^{(qt)^\h  x_{n-1}}_{(q/t)^\h x_{n}}  \Dl(y) \dd{y}
$$
explicitly and obtain, up to a constant factor,
$$
\align
&\det\left(
\matrix
&\vdots&\\
\hdots &
(\xs_i )^{l_j} + (\xs_i )^{-l_j}
-(\xs_{i+1} )^{l_j}- (\xs_{i+1} )^{-l_j}&
\hdots\\
&\vdots&
\endmatrix
\right)_{1\le i,j \le n-1}\,.\\
\intertext{This $(n-1)\times(n-1)$ determinant is equal to
the following $n\times n$ determinant} 
&\det\left(
\matrix
(\xs_1)^{l_1}+(\xs_1)^{-l_1}&\hdots&(\xs_1)^{l_{n-1}}+(\xs_1)^{-l_{n-1}}&1\\
\vdots&&\vdots&\vdots\\
(\xs_n)^{l_1}+(\xs_n)^{-l_1}&\hdots&(\xs_n)^{l_{n-1}}+(\xs_n)^{-l_{n-1}}&1
\endmatrix
\right)
\,.
\endalign
$$
Note that the result is an element of   
$\Lp{-+}[x]$.

Denote by $\I$ the integral
$$
\I=\int_{y\prs x} \b(d_qy) \, \Pm(y; q,t,s\th)\,.
$$
By the above considerations, we have
$$
\I\in\Lp{-+}[x]\,.
$$
Since \tht{A.1} is always divisible by $(u-v)$ it follows
from \tht{3.2} that $\I$ is divisible, for example,  by
$$
(qt\inv \xs_1 - \xs_2) \dots
(q\inv \xs_1 - \xs_2)(\xs_1 - \xs_2)(q\xs_1 - \xs_2) \dots
(tq\inv \xs_1 - \xs_2)\,.
$$
By symmetry, $\I$ is divisible by $D(\xs;q,t)$, and 
moreover
$$
\frac{\I}{D(\xs;q,t)} \in \Lt[x] \,.
$$

Finally, observe that
$$
\alignat2
&\deg V(\ys;q,t)&&=n(n-1)/2\,,\\
&\deg \Pi(\xs;\ys;q,t)&&=2kn(n-1)\,,\\
&\deg D(\xs;q,t)&&=(4k+1)n(n-1)/2 \,.
\endalignat
$$
Therefore, the degree of the RHS of \tht{3.1} is less or
equal to $|\mu|$. This concludes the proof. \qed
\enddemo

\demo{Proof of the theorem 3.1}
Again, we can assume $\theta=2k+1$, $k=0,1,2,\dots$.

One checks the vanishing of the RHS of \tht{3.1}
at the points $x=q^\l$, $\mu\not\subset\l$ by 
precisely the same  argument as in the $A$-series case,
see \cite{Ok3}, \S4. 
If $x=q^\l$ and $\theta$ is an odd natural number, 
then the $q$-integral is just a finite sum all
summand of which vanish. Since the denominator $D(\xs;q,t)$
does not vanish at $x=q^\l$, the equality 
 \tht{3.1} holds up to a constant factor.

To show that this factor equals 1, one can either
compute the RHS at $x=q^\mu$ (there will be only
one non-vanishing summand in the integral), or one
can consider the highest degree term of \tht{3.1}.
This highest term will be
computed explicitly in the next subsection. 
This will conclude the proof of the theorem. \qed
\enddemo

\head 
4.~Highest degree term
\endhead

Denote by
$$
P_\mu(x_1,\dots,x_n;q,t)
$$
the $A$-series Macdonald polynomial with 
parameters $q$ and $t$.

\proclaim{Theorem 4.1 (Highest degree term)} 
$$
\Pm(x_1,\dots,x_n;q,t,s)=
P_\mu(x_1, x_2 t^{-1},\dots,x_n t^{1-n};q,t)+\dots\,, \tag 4.1
$$
where dots stand for lower degree terms.
\endproclaim

By definition, set
$$
\align
&\hV(x_1,\dots,x_n)=\det 
\left[x_i^{n-j}\right]_{1\le i,j \le n}\,, 
\\
&\hPi(x_1,\dots,x_n;y_1,\dots,y_{n-1};q,t) =
\prod_{i,j} \frac {(q y_j/x_i)_\infty}
{(t y_j/x_i)_\infty}\,,
\\
&\hD(x_1,\dots,x_n;q,t)=\prod_{i<j} (x_i-x_j)
\prod_{i\ne j}
\frac {(q x_i/x_j)_\infty}
{(t x_i/x_j)_\infty}\,, \\
&\hC(\mu,n)=\prod_{i<n} B_q(\mu_i+(n-i)\theta,\theta)\,,
\endalign
$$
The hats mean that these  products are related to the top
homogeneous term of $\Pm(x;q,t,s)$. Set also
$$
\align
\hx_i&=x_i t^{1-i}\,, \quad i=1,\dots,n\,,\\
\hy_i&=y_i t^{-i} \,, \quad i=1,\dots,n-1\,.
\endalign
$$

We shall deduce the theorem 4.1 from the following
$q$-integral representation \cite{Ok3} of $P_\mu(x;q,t)$
$$
P_\mu(x;q,t)=\frac 1{\hC(\mu,n)}\frac 1{\hD(x;q,t)}
\int_{y\prec x} P_\mu(y;q,t)\, \hV(y)\, \hPi(x;y;q,t) \d y \,,
$$
where
$$ 
\int_{y\prec x} \d y =
\int_{x_2}^{x_1} \dots \int_{x_n}^{x_{n-1}} \d y_1 \dots \d y_{n-1}
\,.
$$ 

\demo{Proof}
By proposition 2.1 it suffices to consider the case $\mu_n=0$.
In this case we shall use the $q$-integral representation of $\Pm(x;q,t,s)$
and induction on $n$.

By analytic continuation, we can assume that
$$
t=q^\theta, \quad \theta=2k+1, \quad k=0,1,2,\dots\,.
$$

Remark that the highest degree (in variables $u$ and $v$) term of
the polynomial
$$
\prod (1-a u^\1 v^\1)
$$
equals
$$
-a u v (1-a u/v)(1-a v/u)=a^2 u^2 (1- a v /u)(1-a^{-1} v/u)\,.
$$
Therefore
$$
\Pi(\xs;\ys;q,t)= q^{k^2 n(n-1)} 
(t^{n-1} s)^{ 2k n(n-1)} 
\,  \hPi(\hx;\hy;q,t) \prod_i \hx_i^{2k(n-1)}+ \dots
$$
where dots stand for lower degree terms.

Similarly,
$$
V(\ys)= (t^{n-1/2}s)^{n(n-1)/2} \, \hV(\hy) \prod_i \hy_i  + \dots\,,
$$
and
$$
\multline
D(\xs;q,t)=q^{(2k+1)k n(n-1)/2} (t^{n-1} s)^{ (4k+1)n(n-1)/2}
\times\\
 \hD(\hx; q,t) \prod_i \hx_i^{2k(n-1)}+ \dots \,.
\endmultline 
$$
The powers of $q$, $t$, and $s$ in the three above formulas
combine to
$$ 
\frac{
q^{k^2 n(n-1)}(t^{n-1} s)^{ 2k n(n-1)}
(t^{n-1/2}s)^{n(n-1)/2}
}
{  
q^{(2k+1)k n(n-1)/2} (t^{n-1} s)^{ (4k+1)n(n-1)/2}
}
=q^{n(n-1)/4}
$$

Note that
$$
\dd{y_i}= \dd{\hy_i}
$$
and 
$$
\int_{x_{i+1}}^{q x_i} \dd{y_i}\,\big(\dots\big) =
\int_{x_{i+1}}^{t x_i} \dd{y_i}\,\big(\dots\big) =
\int_{\hx_{i+1}}^{\hx_i} \dd{\hy_i}\,\big(\dots\big)\,,
$$
where the first equality is based on \tht{3.2}.

By inductive assumption (note the difference in the definition
of $\hx_i$ and $\hy_i$) we have
$$
\Pm(y_1,\dots,y_{n-1};q,t,s)=t^{|\mu|}
P_\mu(\hy_1,\dots,\hy_{n-1};q,t)
+\dots \,.
$$

Therefore, assuming that the equality \tht{3.1} holds up to
a constant factor $c$, we conclude
$$
\Pm(x_1,\dots,x_n;q,t,s)=
c P_\mu(\hx_1,\dots,\hx_n;q,t)+\dots\,.
$$
But then setting $x_n=1$ and using proposition 2.1,
inductive hypothesis and the stability of $A_n$-series
Macdonald polynomials, we immediately find
$$
c=1\,.
$$
This concludes the proof of the theorem and also the proof
of the theorem 3.1 of the previous section.
\qed
\enddemo

\head
5.~Combinatorial formula
\endhead

Recall that the polynomials $\Pm(x;q,t,s)$ form a linear
basis in $\Lt$.

\definition {Definition 5.1}
Let the polynomials 
$$
\p_{\mu,\nu} (u;n)\in\k[u^\1]
$$
be the coefficients in the following expansion
$$
\Pm(u,x_2,\dots,x_n;q,t,s)=\sum_\nu \p_{\mu,\nu}(u;n) \,
P^*_\nu(x_2,\dots,x_n;q,t,s)\,. \tag 5.1
$$
\enddefinition
Write
$$
\nu\prec\mu
$$
if $\mu_1\ge\nu_1\ge\mu_2\ge\nu_2\ge\dots\ge\nu_{n-1}\ge\mu_n$.

\proclaim{Theorem 5.1 (Branching rule)} We have
$$
\multline
\p_{\mu,\nu}(u;n)=\\ \psi_{\mu/\nu}\, t^{-|\nu|} 
 \prod_{\sq\in\mu/\nu}
\left(u-q^{a'(\sq)} t^{-l'(\sq)} \right)
\left(1-\frac1{ q^{a'(\sq)} t^{2n-2-l'(\sq)} s^2 u } \right)\,, 
\endmultline \tag  5.2 
$$
provided $\nu\prec\mu$ and $\p_{\mu,\nu}(u;n)=0$ otherwise.
Here $\psi_{\mu/\nu}$ are the same weights that appear in
the branching rule for ordinary $A$-type Macdonald polynomials.
\endproclaim

Explicit formulas for the weights $\p_{\mu/\nu}$ are given
in \cite{M}, VI.6.24 and Example VI.6.2.(b); the last formula
is reproduced in \tht{0.4}. We will  not use  explicit formulas
for $\p_{\mu/\nu}$ in the proof of \tht{5.2}.

The branching rule immediately results in the following combinatorial
formula for $\Pm(x;q,t,s)$. Call a tableau $T$ on a diagram $\mu$
a {\it reverse tableau} if its entries strictly decrease
down the columns and weakly decrease in the rows. Denote
by  $T(s)$  the entry  of $T$ in the square $s\in\mu$.

\proclaim{Theorem 5.2 (Combinatorial formula)} We have
$$
\multline
\Pm(x_1,\dots,x_n;q,t,s)= \sum_T \psi_{\T} \prod_{s\in\mu}
t^{1-\T(s)}\times
\\
\left(x_{\T(s)}-q^{a'(\sq)} t^{-l'(\sq)} \right)
\left(1-\frac1{ q^{a'(\sq)} t^{2(n-\T(s))-l'(\sq)} s^2 x_{\T(s)} }
\right)\,,
\endmultline \tag 5.3
$$
where the sum is over all reverse tableaux on $\mu$ with entries
in $\{1,\dots,n\}$.
\endproclaim

Here 
$$
\psi_{\T}\in\C(q,t)
$$
is the same $(q,t)$-weight of a tableau
which enters the combinatorial formula for ordinary Macdonald polynomials
(see \cite{M},\S VI.7)
$$
P_\mu(x;q,t) = \sum_T \psi_{\T} \prod_{s\in\mu} x_{\T(s)} \,.  
$$

Comparing this combinatorial formula with the combinatorial
formula \cite{Ok3} for the 2-parametric interpolation Macdonald polynomials
$\Pm(x;q,t)$ one obtains the following proposition, in which part (b) follows
from \tht{2.2}. 

\proclaim{Proposition 5.3} 
$$
\alignat3
&\text{\rm a)} \qquad &&\Pm(x;q,t,s)\to \Pm(x;q,t)\,,&&\quad s\to\infty\,,\\
&\text{\rm b)} \qquad &
s^{-2|\mu|}&\Pm(x;q,t,s)\to t^{(2-2n)|\mu|}\Pm(1/x;1/q,1/t)\,,&&\quad
s\to 0\,.
\endalignat
$$
\endproclaim

In the proof of the branching rule 
we shall induct on $n$ and use the following
corollary of this theorem

\proclaim{Corollary 5.4 of theorem 5.1} For all $r=1,\dots,n$ we have
$$
\Pm(x_1,\dots,x_r,q^{\mu_{r+1}},\dots,q^{\mu_n};q,t,s) =
c x_1^{\mu_1} \dots x_r^{\mu_r} + \dots\,, \tag 5.4
$$
where $c$ is a non-zero constant  
and dots stand for lower monomials in lexicographic order.
\endproclaim

\demo{Proof} 
Given any partition $\nu$ set
$$
\i\nu=(\nu_i,\nu_{i+1},\dots)\,.
$$
Set also
$$
\pmu=\ii1\mu \,.
$$
Consider the leading term of $\Pm(x;q,t,s)$ as of a polynomial
in $x_1$. Then \tht{5.2} asserts that this leading term
equals  
$$
x_1^{\mu_1}  t^{-|\pmu|} 
P^*_{\pmu} (x_2,\dots,x_n;q,t,s) \,, \tag 5.5
$$
because  $\psi_{\mu/\pmu}=1$. 
Using \tht{5.2} again one obtains the leading term of \tht{5.5}
in $x_2$ and so on. Finally, observe that
$$
P^*_{\i\mu}(q^{\i\mu})\ne 0 
$$
for all $i$. \qed
\enddemo

We shall also need the following elementary

\proclaim{Lemma 5.5} Let $\l$ be a partition. Then
$$
\Pm(x_1,\dots,x_{i-1},q^{\l_i},\dots,q^{\l_n};q,t,s)=0
\quad\text{unless}\quad \i\mu\subset\i\l\,.
$$
\endproclaim
\demo{Proof}
By definition of $\Pm(x;q,t,s)$, this polynomial
vanishes at all points
$$
x_1=q^{\nu_1},\dots,x_{i-1}=q^{\nu_{i-1}}\,,
$$
where 
$$
\nu_1\ge\dots\ge\nu_{i-1}\ge\l_i
$$
are arbitrary integers, provided $\i\mu\not\subset\i\l$. \qed
\enddemo

\demo{Proof of theorem 5.1} Induct on $n$. The case
$n=1$ is clear. 

Fix some $i$. Show that if $\nu_i < \mu_i$  then
$$
\p_{\mu,\nu}(u;n)=0\,, 
\quad u=q^{\nu_i}t^{1-i}, 
\dots, q^{\mu_i-1} t^{1-i}\,. \tag 5.6
$$
We shall prove \tht{5.6} by induction on the partition $\i\nu$,
that is we shall deduce \tht{5.6} from the assumption that 
$$
\p_{\mu,\eta}(u;n)=0\,,
\quad u=q^{\eta_i}t^{1-i}, 
\dots, q^{\mu_i-1} t^{1-i}\,,
$$
for all $\eta$ such that 
$$
\i\eta\subsetneqq\i\nu\,.
$$

Suppose 
$$
u=q^k t^{1-i}, \quad k=\nu_i,\dots,\mu_i-1 \,,
$$
and consider the expansion
$$
\multline
\Pm(q^k t^{1-i},x_2,\dots,x_i,q^{\nu_{i}},\dots,q^{\nu_{n-1}};q,t,s)=\\
=
\sum_{\eta}
\p_{\mu,\eta}(q^k t^{1-i};n)\,  
P^*_\eta(x_2,\dots,x_i,q^{\nu_{i}},\dots,q^{\nu_{n-1}};q,t,s)\,.
\endmultline \tag  5.7
$$
By the lemma  only summands satisfying
$$
\i\eta\subset\i\nu
$$
are nonzero. On the other hand, if
$$
\i\eta\subsetneqq\i\nu
$$
then in particular
$
\eta_i\le\nu_i
$
and by our assumption $\p_{\mu,\eta}(q^k t^{1-i};n)=0$.
Therefore only summands with
$$
\i\eta=\i\nu
$$
enter the sum. 

By the corollary 5.4
applied to polynomial
$$
P^*_\eta(x_2,\dots,x_i,q^{\nu_{i}},\dots,q^{\nu_{n-1}};q,t,s)
$$
in $n-1$ variable,  each summand in \tht{5.7} has the following 
form
$$
c_\eta  \p_{\mu,\eta}(q^k t^{1-i};n)\,
 \big( x_2^{\eta_1}\dots x_i^{\eta_{i-1}} +\dots
\big)\,, \tag  5.8
$$
where $c_\eta$ is a nonzero factor and dots stand for
lower monomials in lexicographic order. 

On the other hand, by the symmetry of $\Pm(x;q,t,s)$ we have
$$
\multline
\Pm(q^k t^{1-i},x_2,\dots,x_i,q^{\nu_{i}},
\dots,q^{\nu_{n-1}};q,t,s)=\\
=\Pm(x_2/t,\dots,x_i/t,q^k,q^{\nu_{i}},\dots,q^{\nu_{n-1}};q,t,s)=0
\endmultline 
$$
by lemma 5.5, provided $\nu_i\le k < \mu_i$. 
Therefore the polynomial \tht{5.7} is identically zero.
By virtue of \tht{5.8} it is impossible unless 
$$
\p_{\mu,\eta}(q^k t^{1-i};n)=0
$$
for all $\eta$ such that $\i\eta=\i\nu$. This proves \tht{5.6}.

Since $\p_{\mu,\nu}(u;n)$ is invariant under the 
transformation
$$
u \mapsto \frac{1}{t^{2n-2} s^2 u}
$$
it vanishes also at the points
$$
\frac{1}{t^{2n-2} s^2 u}=q^{\mu_i-1} t^{1-i}, \dots, q^{\nu_i}
t^{1-i}\,,
$$
for all $i$ such that $\nu_i < \mu_i$.

Now show that 
$$
\p_{\mu,\nu}(u;n)=0
$$
if $\nu\not\subset\mu$. Suppose that 
$$
\nu\not\subset\mu
$$
and $\p_{\mu,\nu}(u;n)$ is not identically zero. Since we know
some zeroes of the polynomial $\p_{\mu,\nu}(u;n)$ we have
$$
\deg \p_{\mu,\nu}(u;n) \ge \sum_{i} \max\{\mu_i-\nu_i,0\}\,.
$$
Hence
$$
\deg  \p_{\mu,\nu}(u;n) \, P^*_\nu(x_2,\dots,x_n;q,t,s) \ge 
\sum_{i} \max\{\mu_i,\nu_i\} > |\mu| \,.
$$
Therefore such a summand cannot occur in the expansion \tht{5.1}.

Thus we can assume that
$$
\nu\subset\mu\,.
$$
Since we know $2|\mu/\nu|$ distinct zeroes of $\p_{\mu,\nu}(u;n)$
and again
$$
\deg \p_{\mu,\nu}(u;n) \le |\mu/\nu|
$$
this polynomial should up to a scalar factor equal
$$
\prod_{\sq\in\mu/\nu}
\left(u-q^{a'(\sq)} t^{-l'(\sq)} \right)
\left(1-\frac1{ q^{a'(\sq)} t^{2n-2-l'(\sq)} s^2 u } \right)\,.
$$
This factor equals $\psi_{\mu/\nu}\, t^{-|\nu|}$ because 
the highest degree term of 
$$
\Pm(x_1,\dots,x_n;q,t,s)
$$
is the $A$-series Macdonald polynomial 
$$
P_\mu (x_1,x_2 t^{-1},\dots,x_n t^{1-n}; q,t)\,.
$$
This concludes the proof. \qed
\enddemo

\head
6.~Pieri formulas and Cauchy identity
\endhead

\definition{Definition 6.1}
Let $\p'_{\l,\mu}(u,n)$ be  the coefficients in the 
following expansion
$$
\multline
\prod_{i=1}^n 
\left(u+ x_i t^{1-i}\right)
\left(1+\frac1{s^2 t^{2n-i-1} u x_i}\right)
\,
\Pm(x;q,t,s)=\\
=
\sum_\l \p'_{\l,\mu}(u,n) \, P^*_\l(x;q,t,s)\,.
\endmultline \tag 6.1
$$
\enddefinition 

It is clear that $\p'_{\l,\mu}(u,n)$ is a polynomial in $u$ and $u^{-1}$
invariant under the transformation
$$
u\mapsto \frac1{s^2 t^{2n-2} u }\,.
$$
Denote by $\mu+\bar1$ the following partition
$$
\mu+\bar1=(\mu_1+1,\dots,\mu_n+1)\,.
$$

Denote by  $\p'_{\l/\mu}$  the coefficients of the Pieri formula
for ordinary $A_n$-type Macdonald polynomials (see \cite{M}, VI.6)
$$
\prod_{i=1}^n 
\left(u+ x_i \right)
\,
P_\mu(x;q,t)=
\sum_\l  \p'_{\l/\mu} \, u^{|\mu+\bar1/\l|} \, P_\l(x;q,t)\,. \tag 6.2
$$
\proclaim{Theorem 6.1 (Pieri formulas)} We have
$$
\multline
\p'_{\l,\mu}(u,n)= \\ \p'_{\l/\mu} 
\prod_{\sq\in \mu+\bar1/\l} 
\left(u+q^{a'(\sq)} t^{-l'(\sq)}\right)
\left(1+\frac1{s^2 q^{a'(\sq)} t^{2(n-1)-l'(\sq)} u} \right)\,, 
\endmultline \tag 6.3
$$
provided
$$
\mu\subset\l\subset\mu+\bar1\,,
$$
and $\p'_{\l,\mu}(u,n)=0$ otherwise.
\endproclaim

Pieri formulas for the interpolation Jack polynomials were considered in \cite{KS}.  

\demo{Proof}
First show that 
$$
\p'_{\l,\mu}(u,n)=0 \quad\text{if}\quad \mu\not\subset\l \tag 6.4 
$$
by induction on $\l$. Assume that 
$$
\p'_{\eta,\mu}(u,n)=0
$$
for all partitions $\eta$ such that
$
\eta \subsetneqq \l\,.
$
Then evaluation of \tht{6.1} at $x=q^\l$ gives \tht{6.4}.

From now on we suppose 
$$
\mu\subset\l\,.
$$
Set $\mu'_0=n$ and suppose that
$
\l'_k < \mu'_{k-1}
$
for some $k\ge 1$. Show that
$$
\p'_{\l,\mu}(-q^{k-1}t^{1-i},n)=0\,, \quad i=\l'_k+1,\dots,\mu'_{k+1}\,.
\tag 6.5
$$
Again, we assume that
$$
\p'_{\eta,\mu}(-q^{k-1}t^{1-i},n)=0\,,
\quad i=\eta'_k+1,\dots,\mu'_{k+1}\,.
$$
for all partitions $\eta$ such that
$$
\eta \subsetneqq \l
$$
and evaluate \tht{6.1} at
$$
x=q^\l, \quad u=-q^{k-1}t^{1-i}\,, 
\quad i=\l'_k+1,\dots,\mu'_{k+1}\,.
$$
This gives \tht{6.5}.

Since $\p'_{\l,\mu}(u,n)$ is invariant 
with respect to the transformation
$$
u\mapsto \frac1{s^2 t^{2n-2} u }\,,
$$
we have
$$
\p'_{\l,\mu}
\left(-\frac1{s^2 q^{k-1}t^{2n-i-1}},n
\right)=0\,, \quad i=\l'_k+1,\dots,\mu'_{k+1}\,. \tag 6.6
$$

Now prove that 
$$
\p'_{\l,\mu}(u,n)=0 \quad\text{if}\quad \l\not\subset\mu+\bar1\,. \tag 6.7 
$$
Suppose that $\p'_{\l,\mu}(u,n)\ne0$.
Then by \tht{6.5} and \tht{6.6}
$$
\deg \p'_{\l,\mu}(u,n) \ge \#\{i\st \mu_i = \l_i\}\,.
$$
If $\l\not\subset\mu+\bar1$ then
$$
\deg \p'_{\l,\mu}(u,n) \,P^*_\l(x;q,t,s) \ge
\sum_i \max\{\mu_i+1,\l_i\} > |\mu|+n \,,
$$
but then such a summand cannot enter the RHS of \tht{6.1}.

If $\l\subset\mu+\bar1$ then by the same reason
$$
\deg \p'_{\l,\mu}(u,n) \le |\mu+\bar1/\l| \,.
$$
By \tht{6.5} and \tht{6.6} we have
$$
\multline \p'_{\l,\mu}(u,n)= \\ \co \, 
\prod_{\sq\in \mu+\bar1/\l} 
\left(u+q^{a'(\sq)} t^{-l'(\sq)}\right)
\left(1+\frac1{s^2 q^{a'(\sq)} t^{2(n-1)-l'(\sq)} u} \right)\,.
\endmultline 
$$
The constant factor is determined by  \tht{4.1} and \tht{6.2}.
This concludes the proof of the theorem. \qed
\enddemo

The similarity of the proof of Pieri formulas 
and of the combinatorial formula is not accidental.
In fact, these two theorems are equivalent to each other
by virtue of the following theorem

\proclaim{Theorem 6.2 (Cauchy identity)}
$$
\multline
\prod_{i=1}^{n} \prod_{j=1}^{m}
\left( x_i t^{n-i} - y_j q^{m-j} \right)
\left( 1-\frac1{s^2 q^{m-j}  t^{n-i} x_i y_j} \right) = \\
=
\sum_{\mu\subset (m^n)}
(-1)^{|\tm|} t^{(n-1)|\mu|} q^{(m-1)|\tm|}\,
\Pm(x;q,t,s)\, P^*_\tm (y;t,q,s) \,, 
\endmultline \tag 6.8
$$
where
$$
\tm=(n-\mu'_m,\dots,n-\mu'_1) \,.
$$
\endproclaim

It follows from  \tht{6.8} that the branching rule  in variables
$x$ is equivalent to Pieri formula in variables $y$ and vice
versa.

\proclaim{Lemma 6.3} Suppose $\mu\subset(m^n)$ and $\l$ has at
most $m$ parts. Then
$$
\prod_{i=1}^n \prod_{j=1}^m
\left( q^{\mu_i} t^{n-i} - t^{\l_j} q^{m-j} \right) = 0 \tag 6.9
$$
unless $\tm\subset\l$.
\endproclaim

\demo{Proof of the lemma} Suppose this product does not vanish.
Then
$$
\l_{m-\mu_i}\ne n-i\,, \quad  i=1,\dots,n\,.
$$
For $i=\mu'_1+1,\dots,n$ we obtain
$$
\l_m \ne 0,1,\dots,n-\mu'_1-1\,.
$$
Therefore
$$
\l_m \ge n-\mu'_1\,.
$$
For $i=\mu'_2+1,\dots,\mu'_1$ we obtain
$$
\l_{m-1} \ne n-\mu'_1,\dots,n-\mu'_2-1\,.
$$
Since $\l_{m-1}\ge\l_{m}\ge n-\mu'_1$ we obtain
$$
\l_{m-1} \ge n-\mu'_2\,.
$$
In the same way we obtain
$$
\l_{m-2} \ge n-\mu'_3,\dots, \l_1 \ge n-\mu'_m\,,
$$
which means $\tm\subset\l$. \qed
\enddemo

\demo{Proof of the theorem}
Denote by $f(x;y)$ the product in the LHS of \tht{6.8}.
Observe that
$$
f(x;y)\in \Lt[x] \quad \text{ and }\quad
f(x;y)\in \L_{q,s}[y] \,.
$$
Expand $f(x;y)$ in polynomials $\Pm(x;q,t,s)$
$$
f(x;y)=\sum_\mu \Pm(x;q,t,s) P^?_\mu(y)
$$
where $P^?_\mu(y)$ are certain elements of $\L_{q,s}[y]$
of degree
$$
\deg P^?_\mu(y) \le |\tm|\,.
$$
Since $m$ is the highest exponent of $x_1$ in $f(x;y)$ only
partitions $\mu$ satisfying
$$
\mu\subset (m^n)
$$
enter the expansion of $f(x;y)$. We have to prove that
$$
P^?_\mu(y) = (-1)^{|\tm|} t^{(n-1)|\mu|} q^{(m-1)|\tm|}\,
P^*_\tm (y;t,q,s) \,. \tag 6.10
$$

Induct on $|\mu|$. Assume that
$$
P^?_\eta(y) = (-1)^{|\teta|} t^{(n-1)|\eta|} q^{(m-1)|\teta|}\,
P^*_\teta (y;t,q,s) \,,
$$
for all partitions $\eta$ such that
$$
\eta\subsetneqq\mu
$$
Then we have
$$
f(q^\mu;y)=\Pm(q^\mu;q,t,s) P^?_\mu(y) +
\sum_{\eta\subsetneqq\mu} P^*_\eta(q^\mu;q,t,s) P^*_\teta(y;t,q,s)\,.
$$
Let $\l$ range over all diagrams with $\le |\tm|$ squares. Then
$$
P^*_\teta(t^\l;t,q,s)=0
$$
because $|\teta|>|\tm|$. By the lemma 
$$
f(q^{\mu};t^\l)=0
$$
unless $\l=\tm$. Therefore
$$
P^?_\mu (t^\l) = 0 
$$
unless $\l=\tm$. This proves \tht{6.10} up to a scalar factor. This
factor is determined by \tht{4.1} and the following formula \tht{6.11} for
ordinary Macdonald polynomials which we recall. \qed
\enddemo

\proclaim{Proposition 6.4} 
$$
\prod_{i=1}^n \prod_{j=1}^m 
(x_i - y_j ) = \sum_{\mu\subset (m^n)} (-1)^{|\tm|}
P_\mu(x;q,t) \,P_\tm(y;t,q) \,. \tag 6.11
$$
\endproclaim
\demo{Proof}
The identity \tht{5.4} in \cite{M}, Ch.~VI, reads
$$
\prod_{i=1}^n \prod_{j=1}^m 
(1+x_i y_j ) = \sum_{\mu\subset (m^n)}
P_\mu(x;q,t) \, P_{\mu'}(y;t,q) \,. 
$$
Therefore \tht{6.11} is equivalent to
$$
P_\tm(y;t,q)= \left(
\prod_i y_i^n \right) P_{\mu'}(1/y;t,q)\,.
$$
It is clear that the RHS of the above equality is a polynomial
in $y$. One easily checks that it is an eigenfunction
of the Macdonald $q$-difference operator $D^1_n$ defined
in \cite{M}, formula \tht{VI.3.4}. \qed
\enddemo

\head 
7.~Binomial formula for Koornwinder polynomials
\endhead

The Koornwinder polynomials \cite{K}
$$
P_\l(x;q,t,a_1,\dots,a_4)\in
\C(q,t,a_1,\dots,a_4)[x_1^\1,\dots,x_n^\1]
$$
are $W$-invariant and depend on
six parameters 
$$
q,t,a_1,a_2,a_3,a_4\,.
$$
These polynomials are orthogonal on the torus
$$
|x_i|=1,\quad i=1,\dots,n
$$
with respect to following measure
$$
\frac 1{(2\pi i)^n} \, \Delta(x)\, \frac{dx}{x} \,,
$$
where 
$$
\Delta(x)= \prod_{i<j} 
\frac{(x_i^\1 x_j^\1)_\infty}{(t x_i^\1 x_j^\1)_\infty}
\prod_i
\frac{(x_i^\1,-x_i^\1,\qh x_i^\1,-\qh x_i^\1)_\infty}
{(a_1 x_i^\1,-a_2x_i^\1,\qh a_3 x_i^\1,-\qh a_4 x_i^\1)_\infty} \,.
$$
Here, by definition,
$$
(u,v,\dots,w)_\infty = (u)_\infty (v)_\infty \dots (w)_\infty 
$$
and
$$
(u)_\infty=(1-u)(1-q u)(1-q^2 u) \cdots \,.
$$
These polynomials specialize to Macdonald polynomials for
classical root systems, see \cite{K} and also \cite{D1}, section 5. 
It is known (\cite{D1}, section 5.2) that the highest degree term
of $P_\l$ is the $A$-type Macdonald polynomial
with parameters $q$ and $t$.

The parameters $a_1,\dots,a_4$ are related to parameters $a,b,c,d$ 
used by Koornwinder via  
$$
a=a_1, \quad b=-a_2, \quad c=\qh a_3, \quad d=-\qh a_4 \,.
$$
It is also convenient to introduce, following J.~F.~van Diejen,
dual parameters $a'_1,\dots,a'_4$ related to parameters
$a_1,\dots,a_4$ by the following  involution
$$
\left(
\matrix
\log a'_1\\ \log a'_2\\ \log a'_3\\ \log a'_4
\endmatrix
\right)
= \frac 12
\left(
\matrix
\format \r&\quad\r&\quad\r&\quad\r \\
1 & 1 & 1 & 1\\
1 & 1 & -1 & -1\\
1 & -1 & 1 & -1\\
1 & -1 & -1 & 1
\endmatrix
\right)
\left(
\matrix
\log a_1\\ \log a_2\\ \log a_3\\ \log a_4
\endmatrix
\right)
\,.
$$
In particular,
$$
a'_1=\sqrt{a_1 a_2 a_3 a_4} \,.
$$
To simplify notation, we shall sometimes omit the
six parameters and write simply  $P_\l(x)$.

A $q$-difference operator whose eigenfunctions are
$P_\l(x;q,t,a_1,\dots,a_4)$ was found by T.~Koornwinder in 
\cite{K}. In \cite{D1} van Diejen found explicit $q$-difference
operators $D_k$, $k=1,\dots,n$ such that
$$
D_k P_\l = E_k(q^\l) P_\l
$$
where
$$
E_k\in\L_{t,a'_1}
$$
and the highest degree term of $E_k$ is 
(up to a constant factor) the $k$-th elementary
symmetric function in variables $q^{\l_i} t^{n-i}$.

The operators $D_1,\dots,D_n$ have the following crucial property.
Define the vectors $\rho$ and $\rho'$ by
$$
\align
q^\rho&=(t^{n-1} a'_1, \dots, t a'_1,a'_1)\,,\\ 
q^{\rho'}&=(t^{n-1} a_1,\dots, t a_1, a_1) \,.
\endalign
$$
Suppose $\l$ is a partition and consider the number
$$
\left[ D_k P_\l \right] (q^{\mu+\rho'})\,. \tag  7.1 
$$
Since $D_k$ is a $q$-difference operator this number is a
combination of values of $P_\l$ at several neighboring  points
with some coefficients which do not depend on $\l$.
In fact, see Lemma 4.3 in \cite{D2}, only points of the form
$$
x=q^{\nu+\rho'}\,,
$$
where $\nu$ is a partition and 
$$
|\nu_i-\mu_i|\le 1,\quad \sum |\nu_i-\mu_i|\le k 
$$
contribute to \tht{7.1}. Moreover, the contribution of
the point
$$
\nu=(\mu_1+1,\dots,\mu_k+1,\mu_{k+1},\dots,\mu_n)
$$
is nonzero.

In the context of Cherednik's double affine Hecke algebra
such a property of $q$-difference operators can be 
proved in an abstract fashion, see \cite{C1,2}.
One can apply the Cherednik's algebra techniques to
Koornwinder polynomials, see \cite{M2,No1}.

This property results in the following theorem

\proclaim{Theorem 7.1 (Binomial formula)}
$$
\multline
\frac{P_\l(x q^{\rho'};q,t,a_1,\dots,a_4)}
{P_\l(q^{\rho'};q,t,a_1,\dots,a_4)} = \\
\sum_{\mu\subset\l} t^{(n-1)|\mu|} a_1^{|\mu|}
\frac{\Pm(q^\l;q,t,a'_1) \Pm(x;q,t,a_1)}
{\Pm(q^\mu;q,t,a'_1) P_\mu(q^{\rho'};q,t,a_1,\dots,a_4) }
\,.  
\endmultline 
$$
\endproclaim

Recall that an explicit product expression for $\Pm(q^\mu;q,t,a'_1)$
is given in definition 1.2 in section 1.

\demo{Proof}
In the same way as in the $A$-series case (see the proof
of the main theorem in \cite{Ok4}) it follows
from the above properties of the operators $D_k$ that
for any partition $\mu$ there exists an $q$-difference
operator $D_\mu$ such that
$$
\left[ D_\mu P_\l \right] (q^{\rho'}) = P_\l(q^{\mu+\rho'}) \,.
$$
This operator $D_\mu$ is a polynomial in $D_1,\dots,D_n$
and 
$$
D_\mu \, P_\l = d_\mu(q^\l) \,P_\l\,,
$$
where 
$$
d_\mu\in\Lambda_{t,a'_1}
$$
and
$$
\deg d_\mu \le |\mu| \,.
$$

Consider the Newton interpolation of the function
$$
f(x)=\frac{P_\l(x q^{\rho'};q,t,a_1,\dots,a_4)}
{P_\l(q^{\rho'};q,t,a_1,\dots,a_4)}
$$
with knots
$$
x=q^\mu\,,
$$
where $\mu$ ranges over partitions with at most $n$ parts.
This Newton interpolation has the form
$$
f(x)=\sum_\mu b_{\mu,\l} \Pm(x;q,t,a_1)\,,
$$
where $b_{\mu,\l}$ is a linear combination of the values
$f(q^\nu)$ for $\nu\subset\mu$. We have
$$
f(q^\nu)=d_\nu(q^\l)\,.
$$
Therefore, 
$$
b_{\mu,\l}=b_\mu(q^\l)
$$
for some polynomial
$$
b_\mu\in\Lambda_{t,a'_1}
$$
of degree
$$
\deg b_\mu \le |\mu| \,.
$$

Since the highest degree term of $f(x)$ equals
$$
\frac{t^{(n-1)|\l|} a_1^{|\l|}}
{P_\l(q^{\rho'};q,t,a_1,\dots,a_4)}
P_\l(x_1,x_2 t^{-1}, \dots , x_n t^{1-n}; q,t)\,,
$$
where $P_\l(x;q,t)$ is the ordinary $A$-type Macdonald
polynomial, and the highest degree term of $\Pm(x;q,t,a_1)$
equals
$$
P_\mu(x_1,x_2 t^{-1}, \dots , x_n t^{1-n}; q,t)
$$
we have
$$
b_\mu(q^\l)=\cases
0, & |\lambda|\le |\mu|, \lambda\ne\mu\,, \\
t^{(n-1)|\mu|} a_1^{|\mu|} /
{P_\mu(q^{\rho'};q,t,a_1,\dots,a_4)}, & \lambda=\mu \,.
\endcases
$$
Since $\deg b_\mu \le |\mu|$ the polynomial $b_\mu$ is
completely determined by its values at the points
$q^\l$, where $|\l|\le|\mu|$.
Therefore $b_\mu(q^\l)$ is proportional to $\Pm(q^\l;q,t,a'_1)$
and precisely 
$$
b_\mu(q^\l) = \frac{t^{(n-1)|\mu|} a_1^{|\mu|}}
{P_\mu(q^{\rho'};q,t,a_1,\dots,a_4)} \frac{\Pm(q^\l;q,t,a'_1)}
{\Pm(q^\mu;q,t,a'_1)}\,.
$$
This concludes the proof. \qed
\enddemo

The following important  property of Koornwinder polynomials 
immediately follows from the binomial theorem

\proclaim{Theorem (\rm J.~P.~van Diejen, \cite{D2})} If $a'_1=a_1$ then 
$$
\frac{P_\l(q^{\nu+\rho'};q,t,a_1,\dots,a_4)}
{P_\l(q^{\rho'};q,t,a_1,\dots,a_4)}=
\frac{P_\nu(q^{\l+\rho};q,t,a_1,\dots,a_4)}
{P_\nu(q^{\rho};q,t,a_1,\dots,a_4)} \,.
$$
\endproclaim

The general symmetry,
$$
\frac{P_\l(q^{\nu+\rho'};q,t,a_1,\dots,a_4)}
{P_\l(q^{\rho'};q,t,a_1,\dots,a_4)}=
\frac{P_\nu(q^{\l+\rho};q,t,a'_1,\dots,a'_4)}
{P_\nu(q^{\rho};q,t,a'_1,\dots,a'_4)} \tag 7.2
$$
conjectured in \cite{D2}, depends on a
formula for $P_\l(q^{\rho'};q,t,a_1,\dots,a_4)$.
The conjectural formula for this number 
(see formula \tht{5.5} in \cite{D2})
was proved in \cite{D2} 
under the self-duality condition
$$
a'_1=a_1\,.
$$

According to the note added in proof to \cite{D2},
that formula (together with symmetry \tht{7.2}) 
was proved recently by I.~G.~Macdonald (in preparation).

\head Appendix 1. $q$-Integrals 
\endhead
Recall that (see, for example, Ch.~1 of \cite{GR})
$$
\int_v^u z^l \dd z = \frac1{[l]_q} (u^{l} - v^{l}),
\quad l\ne 0\,, \quad [l]_q=\frac{1-q^l}{1-q}\,. \tag A.1
$$
If $u=q^k v$ for some $k\in\N$ then
$$
\int_v^u f(z) \dd{z}= - (1-q)  \sum_0^{k-1} f(v q^i) 
$$
for any polynomial $f(z)$ without constant
term. For example,  polynomials 
$f(z)\in\C[z^\1]$ which satisfy
$$
f(z\inv )=-f(z)
$$
do not have constant terms.
Suppose an analytic function 
$f(z)$ is given by its Laurent series of the form
$$
f(z)=\sum_1^\infty c_n(z^n-z^{-n}),\quad \de\inv<|z|<\de, 
$$
for some $\de>1$.
By definition, set
$$
\align
\int_v^u f(z)\dd{z}&=
(1-q) \sum_1^\infty c_n \left.\left(
\frac1{1-q^n}z^n - \frac1{1-q^{-n}}z^{-n} \right) 
\right|_{z=v}^{z=u} \\
&=
(1-q) \sum_1^\infty \left. \frac{c_n}{q^{-n/2}-q^{n/2}} \left(
(z \qm)^n + (z \qm)^{-n} \right) 
\right|_{z=v}^{z=u} \,, \tag  A.2
\endalign
$$
which converges if
$$
|q|\de\inv < |v|,|u| < \de \,.
$$
 
Note the following symmetries of \tht{A.2}
$$
\int_v^u f(z)\dd{z}=
\int_v^{q/u} f(z)\dd{z}=
\int_{q/v}^u f(z)\dd{z}=
-\int_u^v f(z)\dd{z}\,.
$$

Recall that the $q$-beta function is defined (see \cite{GR},1.10) by
$$
B_q(a,b)=(1-q)\frac{(q)_\infty (q^{a+b})_\infty}
{(q^{a})_\infty (q^{b})_\infty}  \,. \tag A.3
$$

\head
Appendix 2. Absence of rational $q$-difference equations 
\endhead     

In this Appendix we show that the following polynomials
$f(x)$ in one variable $x$
$$
\align
f_m(x)&=P^*_m(x;q,t,s)\\
&=(x-1)(x-q)\cdots(x-q^{m-1}) 
\left(1-\frac1{s^2 x}\right)
\cdots  
\left(1-\frac1{s^2 q^{m-1} x}\right)\,,
\endalign
$$
where
$$
m=1,2,\dots\,,
$$
do not satisfy any $q$-difference equation of the
form
$$
\sum_{-d\le i\le d} a_i(x) f_m(q^i x) = E(m) f_m(x)\,, 
\quad d>0\,, \tag A.4
$$
where $a_i(x)$ are rational functions in $x$ and
$$
a_d(x)\not\equiv 0
\quad\text{ or }\quad
a_{-d}(x)\not\equiv 0\,.
$$

Assume the polynomials $f_m(x)$ satisfy \tht{A.4}. Consider the
points
$$
x=q^{m-d}, \quad m=2d,2d+1,\dots \,.
$$
Only finitely many of them are poles of functions
$a_i(x)$, $-d\le i\le d$. Evaluate \tht{A.4} at the 
remaining points. We obtain
$$
a_d(q^{m-d}) f_m(q^m) = 0 \,.
$$
Since $f_m(q^m)\ne 0$ we obtain
$$
a_d(q^{m-d})=0
$$
for infinitely many values of $m$. Therefore
$$
a_d(x)\equiv 0 \,.
$$
Similarly, evaluation at the points 
$$
x=\frac1{s^2 q^{m-d}}, \quad  m=2d,2d+1,\dots
$$
gives
$$
a_{-d} \left(\frac1{s^2 q^{m-d}}\right) = 0
$$
for infinitely many values of $m$. Therefore
$$
a_{-d}(x)\equiv 0 \,.  
$$

\enddocument

\end